\documentclass[11pt]{amsart}
\usepackage[english]{babel}
\usepackage{amssymb}
\usepackage{fourier}
\usepackage{mathrsfs}
\usepackage{enumerate}
\usepackage{color}
\usepackage{ifpdf}
\usepackage{tikz}
\usepackage{tikz-cd}
\usetikzlibrary{decorations.pathmorphing,arrows}
\usepackage[initials]{amsrefs}

\usepackage{caption}
\usepackage{subcaption}
\usepackage{comment}
\usetikzlibrary{intersections}

\newcommand{\bbZ}{{\mathbb Z}}

\newcommand{\SL}{\operatorname{SL}}

\newcommand{\Isom}{\operatorname{Isom}}

\newcommand{\Homeo}{\operatorname{Homeo}}

\newtheorem{theorem}{Theorem}[section]

\newtheorem{lemma}[theorem]{Lemma}

\newtheorem{corollary}[theorem]{Corollary}

\theoremstyle{definition}
\newtheorem{definition}[theorem]{Definition}

\newtheorem{remark}[theorem]{Remark}

\numberwithin{equation}{section}

\author{Yanlong Hao}
\address{University of Michigan, Ann Arbor}
\email{ylhao@umich.edu}


\begin{document}
\date{}
\title{Groups that (do not) act isometrically on hyperbolic spaces}
\maketitle
$\textbf{Abstract}$: In this paper, we show that if a group acts isometrically on a good hyperbolic space of finite volume entropy through a non-elementary action, then it admits an affine action on some $L^p$-space with an unbounded orbit for sufficiently large $p$. As an application, we prove that any isometric action of a group with the fixed point property $F_\infty$ on a good hyperbolic space must have a bounded orbit.

\section{Introduction}
Fixed point properties for group actions on Banach spaces have played a central role in modern geometric group theory and rigidity theory. Beginning with Kazhdan's seminal work on Property (T) \cite{MR209390}, which guarantees fixed points for affine isometric actions on Hilbert spaces, much effort has been devoted to extending these rigidity phenomena to broader classes of Banach spaces, including 
 $L^p$-spaces and uniformly convex spaces. These properties often serve as powerful obstructions to nontrivial group actions on various geometric and analytic objects.

Recent work has revealed significant connections between fixed point properties on Banach spaces and the dynamics of group actions on negatively curved or hyperbolic spaces. For example, Bader, Furman, Gelander, and Monod \cite{MR2316269} developed a framework for understanding group actions on uniformly convex Banach spaces, including $L^P$-spaces, and established fixed point criteria in broad generality. Separately, Monod \cite{MR2715399}*{Chapters 13 and 14} and others have explored the role of bounded cohomology in controlling affine isometric actions, with consequences for rigidity.

It is known that hyperbolic groups themselves do not have fixed point properties for actions on $L^p$-spaces for $p$ sufficiently large. Yu \cite{MR2221161} constructed proper affine isometric actions of hyperbolic groups on $L^p$-spaces for sufficiently large $p$, and Nica \cite{B} later provided a different construction of such actions. Both constructions are based on the work of Mineyev \cites{MR1866802, I}, respectively.

The isometry group of a proper hyperbolic space shares many structural and dynamical properties with hyperbolic groups. In this work, we extend the results of Yu and Nica by showing that similar phenomena hold for unbounded subgroups of the isometry group of proper hyperbolic spaces under mild conditions. 

In this paper, all spaces have finite \textit{volume entropy}. See discussion in Section~\ref{volume growth}.

\begin{theorem}\label{2}
Let $(X,d)$ be a (Gromov) hyperbolic complex or a good (Gromov) hyperbolic space. Assume that $\Isom(X)$ is non-elementary, and $G<\Isom(X)$ is bounded. Then the linear isometric action of $G$ on the space $L^p(\partial X\times \partial X, \nu^{BM})$ has an unbounded cocycle for all $p$ sufficiently large.
\end{theorem}

Here, $\nu^{BM}$ is the Bowen-Margulis measure for the $\Isom(X)$-action; see Section 4 for the construction in this setting. Also, see Section~\ref{section: good space} for the definition of good hyperbolic spaces.

Following \cite{MR2316269}, we say that a topological group has the \textit{fixed point property $F_p$} for some $1\leq p<\infty$ if every affine isometric action of $G$ on any $L^p$-space admits a global fixed point. We say that $G$ has \textit{fixed point property $F_\infty$} if it has fixed point property $F_p$ for all $1\leq p<\infty$. Note that $F_2$ is equivalent to Kazhdan's property (T) \cites{MR578893, MR340464} for locally compact, $\sigma$-compact groups.

An interesting corollary of Theorem~\ref{2} is the following.
\begin{corollary}\label{4}
Let $G$ be a group with fixed point property $F_\infty$, and let $X$ be a good hyperbolic space such that $\Isom(X)$ is not elementary. Then any isometric action of $G$ on $X$ has a bounded orbit.
\end{corollary}

There are a few families of groups known with the fixed point Property $F_\infty$:
\begin{enumerate}
    \item Simple higher (real) rank algebraic groups and their lattices \cite{MR2316269}*{Theorem B}. They even have Lafforgue's strong property (T) \cite{MR4018265}. 
    \item The universal lattice $\SL_n(\bbZ[x_1,x_2,\cdots,x_k])$ for $k\geq 0$ and $n\geq 4$ \cite{MR2794627}*{Theorem 1.2}.  
    \item Higher rank Steinberg groups over commutative rings \cite{oppenheim2023banach}.
    \item Certain acylindrically hyperbolic groups arising as quotients of a family of hyperbolic groups, \cite{MO}*{Corollary 1.2}.
\end{enumerate}

For higher rank lattices, Corollary~\ref{4} has been proved by Haettel \cite{MR4094562}, and also by Bader-Caprace-Furman-Sisto \cite{MR4874176} using different methods. 

In \cite{MO}, the authors construct an acylindrically hyperbolic group with property $F_\infty$. Thus, Corollary~\ref{4} does not hold when $X$ is an arbitrary hyperbolic space. 

\medskip

There are several other directions in which the constructions of Yu and Nica have been generalized. For instance, analogous results have been established for relatively hyperbolic groups in \cite{IF}. A similar statement is also claimed in \cite{CD}*{Theorem 0.1} for groups acting on hyperbolic spaces while preserving a non-collapsing Radon measure.

\medskip

\textbf{About the proof.} The proof follows the idea of Nica \cite{B}. However, since we consider spaces that are not necessarily proper or geodesic, there are technical issues that must be addressed.

First, subgroups of $\Isom(X)$ are not necessarily locally compact, so even the construction of the Patterson–Sullivan measure is not straightforward. Here, we utilize the idea that many geometric arguments can be carried out by working with the Gromov boundary, which is compact due to the finiteness of the volume entropy. Note that the Patterson-Sullivan measures have been constructed in \cite{MR3558533}*{Chapter 16} for groups of generalized divergence type.

Secondly, in \cite{B}, the properness of the cocycle is proved under the assumption that the Patterson–Sullivan measure is Ahlfors regular. The Ahlfors regularity in \cite{B} is established by showing that the measure lies in the same class as the Hausdorff measure. In more general settings, a similar property, good enough for the proof, of the Patterson–Sullivan measure is derived using the shadow lemma. Since the essential structure of the proof remains unchanged, we only highlight the necessary modifications required in our setting.

\subsection*{Acknowledgements}
I would like to thank Alexander Furman for his generous guidance and support during the completion of this work.

\section{Preliminaries}
\subsection{Hyperbolic spaces}
As a general reference, we refer the reader to \cite{MA}.  

Let $(X,d)$ be a metric space and $o\in X$ be a fixed base point. Define the Gromov product $(\cdot,\cdot)_o:X\times X\rightarrow[0,\infty)$  as
$$(x,y)_o=\frac{1}{2}(d(x,o) +d(y,o)-d(x,y)).$$ We simply denote $(\cdot,\cdot)_o$ by $(\cdot,\cdot)$.
A different base point $p$ leads to another Gromov product such that
$$|(x,y)_o-(x,y)_p|\leq d(o,p).$$
If there exists $\delta\geq 0$ such that the Gromov product satisfies
\begin{equation}\label{1}
(x,y)\geq \min\{(x,z),(y,z)\}-\delta
\end{equation}
for some (equivalently, for every) base point $o\in X$ and any $x$, $y$, $z\in X$, the space $X$ is said to be \textit{($\delta$-)hyperbolic}.

 A sequence $(x_n)$ in $X$\textit{ tends to $\infty$} if
$$\lim_{i,j\rightarrow \infty}(x_i,x_j)=\infty.$$
Two such sequences $(x_n)$ and $(y_n)$ are equivalent if $\lim_{n\rightarrow\infty}(x_n,y_n)=\infty.$ \textit{The Gromov boundary} of $X$, denoted $\partial X$, is the set of equivalence classes of sequences tending to infinity. The space $\overline{X}=X\cup \partial X$ can be given a natural topology making it a compactification of $X$, on which the isometry group $\textbf{Isom}(X)$ acts by homeomorphisms. 

The Gromov product can be extended (not necessarily continuously) to $\overline{X}$ in such a way that the estimate \eqref{1} holds with different constants.  

For $L$, $C>0$, a map $f:(X,d)\to (Y,d')$ is called a \textit{$(L,C)$-quasi-isometry} if the following two conditions hold:
\begin{enumerate}
    \item 
    $\frac{1}{L}d'(f(x),f(y))-C\leq d(x,y)\leq Ld'(f(x),f(y))+C$
    for all $x$, $y\in X$.
    \item For any $x'\in Y$, there exists $x\in X$ such that $d'(f(x),x')\leq C$.
\end{enumerate}
 A \textit{rough-isometry} is an $(1,C)$-quasi-isometry.
\subsection{Volume entropy of metric spaces}\label{volume growth}
Let $(X,d)$ be a metric space with base point $o$. A subset $E\subset X$ is called \textit{$k$-separated} if $d(x,y)\geq k$ for all $x\neq y\in E$. A $k$-separated set is called \textit{maximal} if $E\cup \{x\}$ is not $k$-separated for all $x\notin E$. For $R\geq 0$, let 
\[
N_R^k=\max_{E}\#E,
\]
where $E$ runs over all $k$-separated subset of $B(o,R)$. The \textit{volume entropy at scale $k$} of $X$ is defined as
\[
\delta_X^k=\limsup_R \frac{\ln N^k_R}{R}.
\]
The \textit{volume entropy} of $X$ is given by 
\[
\delta(X)=\inf_{k>0} \delta_X^k
\]

\begin{lemma}
    If $X$ is hyperbolic and $\delta_X< \infty$, $\partial X$ is compact.
\end{lemma}
\begin{proof}
    By definition, there exists $k>0$ such that for any $E\subset X$ which is a maximal $k$-separated set,
    \[
    \limsup_R \frac{\ln \# E\cap B(o, R)}{R}\leq \delta_X+1.
    \]
 Now $E$ with the induced metric is roughly isometric to the space $X$. Then they share the same Gromov boundary. The result follows from the fact that $E$ is a proper space.
\end{proof}

\subsection{Limit set and non-elementary subgroup}
Let $G$ be a subgroup of the isometry group $\textbf{Isom}(X)$. Consider the set $\Lambda_G$ of accumulation points of the orbit of $o$ in $\partial X$. $\Lambda_G$ is in fact independent of the choice of base point $o$, and we call this set the \textit{limit set} of $G$. The group $G$ is called \textit{non-elementary} if $\Lambda_G$ contains at least three elements. A consequence of $G$ being non-elementary is that it has hyperbolic elements.
\subsection{Good hyperbolic spaces}\label{section: good space} The Gromov product of a hyperbolic space $X$ can be extended to $\bar{X}$. However, the extension is, in general, not a continuous map. To address this,  good hyperbolic spaces have been introduced in \cite{MR3551185}.

\begin{definition}\cite{MR3551185}*{Definition 3.1}
    We say that a hyperbolic space X is $\epsilon$-good, where $\epsilon>0$, if the following two properties hold for each basepoint $o\in X$:
    \begin{enumerate}
        \item the Gromov product $(\cdot,\cdot)$ extends continuously to the bordification $\bar{X}$,
 and
 \item $e^{-\epsilon(\cdot,\cdot)}$ is a metric on the boundary $\partial X$.
    \end{enumerate}
\end{definition}
A space X is called \textit{good hyperbolic} if it is $\epsilon$-good for some $\epsilon$. In this case, the metric $d_\partial=e^{-\epsilon(\cdot,\cdot)}$ is called the \textit{visual metric}.

CAT$(-1)$ space and Green metric from a finitely supported random walk on a hyperbolic group are examples of good hyperbolic spaces. 

Another family of good hyperbolic spaces is the hat metrics $\hat{d}$ constructed by Mineyev \cite{I}. We recall the key property we require in this paper. 

Let $(X,d)$ be a hyperbolic complex with isometry group $\textbf{Isom}(X)$. There is a metric $\hat{d}$ on $X$ with the following properties:
\begin{enumerate}
\item
$\hat{d}$ is $\textbf{Isom}(X)$-invariant.
\item
$(X,d)$ and $(X,\hat{d})$ are quasi-isometric. Hence, the two metrics have homeomorphic boundary $\partial X$.
\item
Fix $p\in X$. The $\hat{d}$-Gromov product $\widehat{(x,y)}=\frac{1}{2}(\hat{d}(p,x)+\hat{d}(p,y)-\hat{d}(x,y))$ on $X\times X$ extends continuously to a function $\overline{X}\times \overline{X}\rightarrow [0,+\infty]$.
\item
There exist $\varepsilon_0>0$ such that $d_\partial(x,y)=e^{-\varepsilon_0\widehat{(x,y)}}$ defines a metric on $\partial X.$
\item
$\textbf{Isom}(X)$ acts on $\partial X$ by M\"{o}bius homeomorphisms.
\end{enumerate}

\subsection{Shadow Lemma}
The Shadow Lemma describes the property of conformal measures. See \cite{MR3558533}*{Chapter 14} for detailed discussions.
\begin{definition}\cite{MR3558533}*{Definition 4.5.1}
    Let $X$ be a hyperbolic space and $o$ a base point. For each $\sigma>0$ and $x\in X$, let the shadow set be
    \[
    \operatorname{Shad}(x,\sigma)=\{\eta\in\partial X\ |\ (o, \eta)_x\leq \sigma\}.
    \]
\end{definition}

\begin{definition}\cite{MR3558533}*{Definition 15.1.1}
    Let $X$ be a hyperbolic space, and $G<\Isom(X)$. For each $s\geq 0$, a nonzero measure $\mu$ on $\partial X$ is called
 $s$-conformal if
 \[\mu(g(A))=\int_A |g'|^s(\xi) d\mu(\xi)\]
\end{definition}
See Lemma~\ref{Lemma: geometric mean} for definition and properties of $|g'|$.

\begin{lemma}\cite{MR3558533}*{Lemma 15.4.1}Let $X$ be a hyperbolic space, and $G<\Isom(X)$. 
    Let $\mu$ be an $s$-conformal measure on $\partial X$ which is not
 a point mass. Then there exist $\sigma_0$, $C_1$ and $C_2$ such that for all $\sigma>\sigma_0$,   $g\in G$,

 \[
 C_1e^{-sd(go,o)}\leq\mu(\operatorname{Shad}(go,\sigma))\leq C_2e^{-sd(go,o)}
 \]
\end{lemma}

Let $G<\Isom(X)$ be a non-elementary subgroup. In Section~\ref{sec: PS measure}, we will construct the Patterson-Sullivan measure $\nu_G$ of $G$. These are all $\delta_G$-conformal measures. 

\subsection*{2.F. Ultrafilters and Ultralimits}

We recall the notion of ultrafilters and ultralimits, which will be used later (in Section~\ref{sec: PS measure}) to define limits of measures when standard compactness arguments are not available due to a lack of properness.

An \emph{ultrafilter} $\mathcal{U}$ on $\mathbb{N}$ is a collection of subsets of $\mathbb{N}$ satisfying the following properties:

\begin{enumerate}
    \item $\emptyset \notin \mathcal{U}$;
    \item If $A \in \mathcal{U}$ and $A \subset B$, then $B \in \mathcal{U}$;
    \item If $A, B \in \mathcal{U}$, then $A \cap B \in \mathcal{U}$;
    \item For every subset $A \subset \mathbb{N}$, either $A \in \mathcal{U}$ or $\mathbb{N} \setminus A \in \mathcal{U}$.
\end{enumerate}

An ultrafilter is called \emph{non-principal} if it contains no finite set. The existence of non-principal ultrafilters follows from Zorn’s Lemma.

Given a topological space $Y$ and a sequence $(x_n)_{n \in \mathbb{N}}$ in $Y$, we say that $x \in Y$ is the \emph{ultralimit} of $(x_n)$ along a non-principal ultrafilter $\mathcal{U}$, denoted
\[
x = \lim_{n \to \mathcal{U}} x_n,
\]
if for every open neighborhood $V$ of $x$, the set $\{n \in \mathbb{N} \mid x_n \in V\}$ lies in $\mathcal{U}$. That is, the sequence eventually stays in every neighborhood of $x$ with respect to the ultrafilter.

If $Y$ is compact and Hausdorff, then such an ultralimit always exists and is unique. In particular, any bounded sequence in a compact metric space admits an ultralimit with respect to any non-principal ultrafilter \cite{MR1253544}*{2.A and the discussion below}.

\section{Coformal structure of $\partial X$}
In this part, we review basics for M\"{o}bius homeomorphisms. The crucial point is that $|g'|$ is Lipschitz for any M\"{o}bius homeomorphism $g$. 

\subsection{M\"{o}bius homeomorphisms}
Throughout, let $(K, d)$ be a compact metric space without isolated points. Recall that the \textit{cross-ratio} of a quadruple of distinct points in $K$ is defined by the formula
$$(z_1,z_2;z_3,z_4)=\frac{d(z_1, z_3)d(z_2, z_4)}{d(z_1, z_4)d(z_2, z_3)}.$$
A homeomorphism $g:K\rightarrow K$ is called a \textit{M\"{o}bius homeomorphism} if $g$ preserves the cross-ratios, i.e., $(gz_1,gz_2;gz_3,gz_4)=(z_1,z_2;z_3,z_4)$ for all quadruples of distinct points $z_1$, $z_2$, $z_3$, $z_4\in K$.

For M\"{o}bius homeomorphisms, we have the following,
\begin{lemma}\cite{B}\label{Lemma: geometric mean}
Let $g$ be a self-homeomorphism of a compact metric space $K$. Then $g$ is M\"{o}bius if and only if there exists a positive continuous function on $K$, denoted $|g'|$, with the property that for all $x, y\in K$ we have
$$d^2(gx, gy)=|g'|(x)|g'|(y)d^2(x, y).$$
\end{lemma}
In fact, for any $x\in K$, we have
$$|g'|(x)=\lim_{y\rightarrow x}\frac{d(gx,gy)}{d(x,y)}.$$

The next lemma shows that metric derivatives are more than just continuous:
\begin{lemma}\cite{B}
Let $g$ be a M\"{o}bius self-homeomorphism of $K$. Then $|g'|$ is Lipschitz.
\end{lemma}
Since $K$ is compact, there exist constants $c>0$, $C>0$ depending on $g$ with $c<|g'|<C$. It implies that $\ln{|g'|}$ is Lipschitz.

\subsection{$\Isom(X)$ as  M\"{o}bius maps}
When $X$ is a good hyperbolic space, $\Isom(X)$ acts on its Gromov boundary by homemorphisms. A routine check of Gromov product shows that this action is via M\"{o}bius self-homeomorphisms with respect to the visual metric $d_\partial$.

Direct calculation shows that in this case 
\begin{equation}\label{eq: derivate to busmann}
|g'|(\xi)=e^{-\varepsilon_0(2(g \xi, g o)-d(o,g o))},
\end{equation}for $g\in \textbf{Isom}(X)$, $\xi\in \partial X$. The cocycle  $\ln|g'|$ up to scale is called the \textit{Busemann cocycle}. More precisely, the Busemann cocycle $b: \Isom(X)\times \partial X\to \mathbb{R}$ is given by 
\[
b(g,\xi)=-\frac{1}{\epsilon_0 }\ln g'(\xi)=2(g\xi,go)-d(o,go).
\]
The following lemma is true.
\begin{lemma}\label{k}
For any $g\in \textbf{Isom}(X)$, there exist $\kappa>0$ depending on $g$ such that, for any $\xi$, $\eta\in \partial X$,
$$|(g \xi, g o)_p-(g \eta, g o)|<\kappa d_\partial(\xi,\eta).$$
\end{lemma}
\section{Patterson-Sullivan measures}\label{sec: PS measure}
Let $X$ be a good hyperbolic space $X$ and $G$ a topological group. Let $\pi: G\to\Isom(X)$ be a non-elementary action. In this section, we construct the \textit{Bowen-Margulis measure} for $G$.  It starts with the construction of the \textit{Patterson-Sullivan measure}. The method used here is similar to the one in \cite{Q}*{Section 4.4}. 

We start with a few lemmas. For readability, we include a few proofs.
\begin{lemma}\cite{hao2024arithmeticity}*{Lemme 3.1}\label{label: kernal has bound orbit} 
    Let $X$ be a good hyperbolic space with $\partial X$ having more than two points. Denote the map $\partial:\Isom(X)\to \Homeo(\partial X)$. Then $\ker(\partial)o$ is bounded.
\end{lemma}
\begin{proof}
    Let $f\in \ker(\partial)$ and $\xi_i$, $1\leq i\leq 3$, be three different points in $\partial X$. Since $f$ is a M\"{o}bius map, by Lemma~\ref{Lemma: geometric mean}, $f'(\xi_i)f'(\xi_{i+1})=1$ for all $1\leq i\leq 3$. It follows that $f'(\xi_i)=1$, $1\leq i\leq 3$.

    Now by equation~\eqref{eq: derivate to busmann}, 
    \[
    2(\xi_i, f(o))=d(o,f(o)).
    \]
    The hyperbolic inequality
    \[
    (\xi_1,\xi_2)\geq \min\{(\xi_1,f(o)), (\xi_2,f(o))\}-\delta
    \]
    implies that $d(o,f(o))\leq 2(\xi_1,\xi_2)+2\delta$.
    
    This completes the proof.
\end{proof}
    Since $\ker(\partial)$ is a normal subgroup of $\Isom(X)$. The same bound holds for $\Isom(X)o$. By replacing $X$ by $\Isom(X)o/\ker(\partial)$, which is roughly isometric to $\Isom(X)o$, we assume from now on that $\ker{\partial}$ is trivial.

\begin{lemma}\cite{B}*{Proposition 5.6}
    Let $\operatorname{M\"{o}b}(K)$ be the group of M\"{o}bius maps of a compact metric space $K$ with the uniform convergence topology. $\operatorname{M\"{o}b}(K)$ is locally compact and $\sigma$-compact.
\end{lemma}
\begin{lemma}
    Let $X$ be a good Gromov hyperbolic space of finite volume entropy. Let $\mathcal{F}(\xi)=\min_{\xi'\in\partial X}\{(\xi,\xi')\}$. Then there exists a constant $M$ such that $\mathcal{F}\leq M$.
\end{lemma}
\begin{proof}
    This follows from two facts: 
    \begin{enumerate}
        \item $\partial X$ is compact.
        \item $<\cdot,\cdot>$ is continuous.
    \end{enumerate}
\end{proof}
\begin{lemma}\label{lemma: bounded}
    Let $X$ be a good Gromov hyperbolic space of finite volume entropy. Then $\partial(\Isom(X))$ is a closed subgroup of $\operatorname{M\"{o}b}(\partial X)$.  
\end{lemma}
\begin{proof}
    Let $\partial f_i\to g$ be a Cauchy sequence in $\operatorname{M\"{o}b}(\partial X)$. Then the metric derivatives $|f_i'(\xi)|$ are uniformly bounded for $\xi\in\partial X$. By equation~\eqref{eq: derivate to busmann}, 
    \begin{equation}\label{eq: bounded}
        2(f_i(\xi),f_i(o))-d(o,f_i(o))
    \end{equation}  
    is uniformly bounded for $i\geq 1$ and $\xi\in \partial X$.

    Now, assume $f_i(o)$ is not bounded. Then there exists a subsequence $f_{j_k}$ such that $f_{j_k}(o)\to \eta\in \partial X$. Without loss of generality, we assume $f_i(o)\to \eta$.

    By Lemma~\ref{lemma: bounded}, there exist $\eta'$ such that $(\eta,\eta')\leq M$. The hyperbolic inequality \cite{MA}*{pp 433}
    \[
    (\eta, \eta')\geq \min\{(\eta,f_i(o)), (\eta',f_i(o))\}-2\delta
    \]
    together with 
    \[
    \lim_i (f_i(o),\eta)=\infty
    \]
    imply that 
    \[
    (\eta', f_i(o))\leq M+2\delta
    \]
when $i$ is sufficiently large.

    Hence, for such $i$,
    \[
    2(f_i(f_i^{-1}(\eta')),f_i(o))-d(o,f_i(o))
    \]
    are unbounded, contradicting equation~\eqref{eq: bounded}. 

    Therefore, $f_i(o)$ is bounded. It follows that there exists a subsequence such that $f_{j'_k}\to f\in\Isom(X)$. By continuity of $\partial$, $\partial f=g$.

    The result follows.
\end{proof}
 It follows that the closure $L$ of the image $\partial\circ\pi(G)$ in $\partial(\Isom(X))$ is locally compact. Denote $\mu$ the Haar measure of $L$ and $\bar{L}$ the pre-image of $L$ in $\Isom(X)$. The construction of Patterson-Sullivan measures can be generalized to $L$.

Define the Poincaré series by
\[
P(s)=\int_{L} e^{-sd(o,\bar{g}o)}d\mu(g)
\]
and the critical exponent of $G$ by 
\[
\delta_L=\inf\{s|P(s)<\infty\}.
\]

\begin{lemma}
    With the notations as above, $\delta_L\leq \delta_X<\infty$.
\end{lemma}
\begin{proof}
     For any $k>0$, let $E^L_k=\{\bar{g}_i o\}_{i=1}^\infty$ be a maximal $k$-separated set of $\bar{L}o$.  Let $\Omega_s=\{\bar{g}\in \bar{L}\ |\ d(o,\bar{g}o)< s\}$. Since $E^L_k$ is maximal, $\bar{g}_i\Omega_{k}$ is a covering of $\bar{L}o$. Since $\mu(\Omega_{k})>0$ as $\Omega_{k}$ contains a neighborhood of the unit element in $L$ for $k$ sufficiently large, we have 
    \[
    \mu(\partial\{\bar{g}\in \bar{L}\ |\ d(o,\bar{g}o)\leq k\})\leq \mu(\Omega_{k})\#(E^L_k\cap B(o, k)).
    \]
    Therefore
    \[
    \delta_L\leq \delta_X
    \]
    from the fact
    \[
    \delta_L=\limsup_{k\to \infty }\frac{\ln (\mu(\partial\{\bar{g}\in \bar{L}\ |\ d(o,\bar{g}o)\leq k\}))}{k}
    \]

\end{proof}
Similar to \cite{Q}*{Lemma 4.9}, for any $n$, define a probability measure $\nu_n^L$ on $\overline{X}$ by
$$\nu_n^L=\frac{1}{\Phi(\delta_L+\frac{1}{n})}\int_{L} e^{-(\delta_L+\frac{1}{n})d(o,\bar{g}o)}h(d(o,\bar{g}o))D_{\bar{g}o} d\mu(g),$$
where $\Phi(t)=\int_{L} e^{-td(o,\bar{g}o)}h(d(o,\bar{g}o))d\mu(g)$; $D_{\bar{g}o}$ is the Dirac measure at $\bar{g}o$; and here $h=1$ when $\Phi(\delta_L)=\infty$.

Next, we aim to find a weakly star limit of $\nu_n^L$. Note that we do not assume $X$ is proper. Therefore, the existence of weakly star limit is not straightforward. Instead, we use the following construction.

Let $\phi\in C(\partial X)$ be a continuous function. We say that $\bar{\phi}$ is a bounded extension of $\phi$ if 
\begin{enumerate}
    \item $\bar{\phi}\in C(\bar{X})$, and $\bar{\phi}|_{\partial X}=\phi$;
    \item $||\bar{\phi}||_\infty=||\phi||_\infty$.
\end{enumerate}

Fix a non-principal ultrafilter $\mathcal{U}$ of $\mathbb{N}$. For any $\phi\in C(\partial X)$, choose a bounded extension $\bar{\phi}$. Since $|\nu_n^L(\bar{\phi})|\leq ||\phi||_\infty$, 
\[
\lim_{n\to \mathcal{U}} \nu_n^L(\bar{\phi})
\]
is well-defined.

Let $\bar{\phi'}$ be another bounded extension of $\phi$. For any $\epsilon$, there exists a neighbourhood $V$ of $\partial X$ so that $|\bar{\phi}-\bar{\phi'}|$ restricted on $V$ is less than $\epsilon$. Since the space $\bar{X}\setminus V$ is bounded with respect to the distance $d$ and $\Phi(\delta_L)=\infty$, for $n$ sufficiently large, $\nu_n^L(\bar{X}\setminus V)\leq \epsilon$. Therefore for $n$ sufficiently large, 
\[
|v_n(\bar{\phi}-\bar{\phi'})|\leq \epsilon(2||\phi||_\infty+1).
\]
Thus, \[
\lim_{n\to \mathcal{U}} \nu_n^L(\bar{\phi})
\]
is independent of the choice of bounded extension of $\phi$. This defines a measure $\nu_L$ on $\partial X$ by 
\[
\nu_L(\phi)=\lim_{n\to \mathcal{U}} \nu_n^L(\bar{\phi}).
\]
We call $\nu_L$ the \textit{$\mathcal{U}$-weak star limit} of $\nu_n^L$.  

By definition, it is easy to verify that $\nu_L(1)=1$. Hence, the limit process has no escape of mass.

It is trivial to see that for any $g\in L$, and $\xi\in \partial X$
$$\frac{g_*d\nu_L}{d\nu_L}(\xi)=e^{\delta_L(2(\xi,g^{-1}o)-d(p,g^{-1}o))}.$$
Hence, it is a Patterson-Sullivan measure of $L$. And, consequently, it is a Patterson-Sullivan measure of $G$.

Now, the Bowen-Margulis measure $e^{2\delta_L(\xi,\eta)}d\nu_L(\xi)d\nu_L(\eta)$ is $L$-invariant, which we will denote by $\nu^{BM}_{L}$. We denote $\nu_{\Isom(X)}^{BM}$ by $\nu^{BM}$.

\section{The cocycle $\beta$}
\subsection{The cocycle is in $L^p$.}\label{Sec: in L^p}
The cocycle is given by
$$\beta_g(\xi, \eta):=(go,\xi)-(go,\eta),$$
for $g\in G$, and $\xi$, $\eta\in \partial X$.

Here, the cocycle is a generalization of the cocycle used in \cite{B}, where $\beta$ is shown to be related to the Busemann cocycle.

Recall that $(go,\cdot)$ is a Lipschitz map on $(\partial X, d)$ for any $g\in G$. That is,
$$|(g \xi, g o)-(g \eta, g o)|<\kappa_g d_\partial(\xi,\eta)=\kappa_g e^{-\varepsilon_0(\xi,\eta)}.$$ 

Let $\tau=\delta(\Isom(X))$. It is clear that the function $|\beta_g|^{\frac{2\tau}{\varepsilon_0}}e^{2\tau(\xi,\eta)}<\kappa_g^{\frac{2\tau}{\varepsilon_0}}$.

For any $p\geq \frac{2\tau}{\varepsilon_0},$ one has
\begin{equation}
\begin{aligned}
||\beta_g||_p^p&=\int_{\partial X\times \partial X}|\beta_g|^{p}d\nu^{BM}\\
&\leq \int_{\partial X\times \partial X} |\beta_g|^{p-\frac{2\tau}{\varepsilon_0}}\kappa_g^{\frac{2\tau}{\varepsilon_0}}d\nu_{\Isom(X)} d\nu_{\Isom(X)}\\
&\leq ||\beta_g||_{\infty}^{p-\frac{2\tau}{\varepsilon_0}}\kappa_g^{\frac{2\tau}{\varepsilon_0}}<\infty
\end{aligned}
\end{equation}

Hence $\beta$ is a cocycle in $L^p(\partial X\times \partial X, \nu_L^{BM})$ for all $p\geq \frac{2\tau}{\varepsilon_0}.$
\subsection{$\beta$ is unbounded}\label{sec: beta unbounded}
The proof of the unboundedness of $\beta$ is similar to the proof of \cite{B}*{Proposition 7.1}. There is one change needed:

 \textit{The claim in the second paragraph for sufficiently large $R$ is given by the Shadow Lemma. And the general case follows from the compactness of $\partial X$.}

 Hence, we have:
 \begin{lemma}\label{lemma: unbound}
     Under the condition of Theorem~\ref{2}, there exist $a$, $b>0$ such that for 
     \[||\beta(g)||_p^p\geq ad(go,o)-b.\]
 \end{lemma}
 \begin{remark}
    If $\Isom(X)$ contains a hyperbolic element and is elementary. Then the limit set of $\Isom(X)$ contains two point $\xi$, $\eta$. Then Dirac measure supported at $(\xi,\eta)$ serves as a Bowen-Margulis measure for $\Isom(X)$. An easy computation shows that in this case, the conclusion of Lemma~\ref{lemma: unbound} is true.
\end{remark}
\subsection{Proof of Theorem~\ref{2}.}
\begin{proof}
    If $(X,d)$ is a hyperbolic complex, we replace it with the hat-metric $\hat{d}$. Hence, we will assume that $X$ is good hyperbolic. Then the result follows from subsections~\ref{Sec: in L^p} and \ref{sec: beta unbounded}.
\end{proof}

\end{document}